# The Game of Band or Bump

**Bruce Levin**                                                                             BL6@COLUMBIA.EDU
*Department of Biostatistics*
*Mailman School of Public Health*
*Columbia University*
*New York, NY 10032, USA*

## Abstract

In this report we generalize the game of Book or Band introduced in Levin (2024) to an arbitrary playing deck with *m* ranks and *s* cards in each rank, for a total of *t=ms* cards. Two events—a *band* or a *bump*—are defined in terms of given integers *l*, *u*, and *s* with $0 \le l \le u \le s$, not necessarily with $l+u=s$. We derive expressions for the joint stopping time distribution and outcome (band or bump) in terms of rectangular event probabilities for central multiple hypergeometric random variables.

**1. Introduction and joint stopping time and outcome distributions.**

Suppose a generalized deck of playing cards has *m* ranks (or *types*) labeled 1,…,*m*, with *s* cards in each rank, for a total of $t=ms$ cards. Let $X^{(n)} = (X_1^{(n)},...,X_m^{(n)})$ denote the tallies of the *m* ranks after *n* cards are drawn. Given integers $0 \le l \le u \le s$, not necessarily with $l+u=s$, we define a *band* as the event $[l \le X_j^{(n)} \le u, \ j=1,...,m]$ and we define a *bump* as the event $[\max_j \{X_j^{(n)}\} \ge u+1]$ (which obviously can occur only if $u<s$). Here we consider what we call the "Band or Bump" game, where cards are drawn sequentially from the well-shuffled deck until either a band or a bump occurs for the first time. We derive expressions for the joint stopping time and outcome (band or bump) distribution in terms of *rectangular event probabilities* for central multiple hypergeometric random variables. The new game specializes to the original game of Book or Band introduced in Levin (2024) when *l*=1 and $u = s-1$ (with *m*=13 and *s*=4), wherein a bump was a called a "book" as it corresponds to having drawn all four cards in some rank. Levin (2024) shows that if a player wins $2 when the game ends with a book but loses $3 when the game ends in a band, then the game is approximately fair, with an expected gain of just under a nickel per game (with *P*[band] ≈ 0.390753). Another pleasant version of Band or Bump employs *m*=4, *s*=13, *l*=5, *u*=8, which can be played with a standard deck of cards using suits as the *m*=4 types. In this version, though, a win of $2 for a *band* and a loss of $3 for a *bump* is close to fair, with an average gain of about three cents per game (with *P*[band] ≈ 0.605984). See Table 1. Equations (5) and (10) below furnish these and other results.

Table 1

Stopping time and outcome distributions for the game of Band or Bump

with $s=13$, $l=5$, $u=8$, and $m=4$

| $n$ | $P[N=n$, band] | $P[N=n$, bump] | $P[N=n]$ | $P[N=n\mid$band] | $P[N=n\mid$bump] |
|---|---|---|---|---|---|
| 9 |  | 0.000000777369 | 0.000000777369 |  | 0.00000197294 |
| 10 |  | 0.00000634550 | 0.00000634550 |  | 0.0000161047 |
| 11 |  | 0.0000287058 | 0.0000287058 |  | 0.0000728546 |
| 12 |  | 0.0000949860 | 0.0000949860 |  | 0.000241072 |
| 13 |  | 0.000256462 | 0.000256462 |  | 0.000650894 |
| 14 |  | 0.000598412 | 0.000598412 |  | 0.00151875 |
| 15 |  | 0.00124932 | 0.00124932 |  | 0.00317073 |
| 16 |  | 0.00238769 | 0.00238769 |  | 0.00605988 |
| 17 |  | 0.00424478 | 0.00424478 |  | 0.0107731 |
| 18 |  | 0.00710151 | 0.00710151 |  | 0.0180234 |
| 19 |  | 0.0112780 | 0.0112780 |  | 0.0286232 |
| 20 | 0.0217752 | 0.0171131 | 0.0388883 | 0.0359336 | 0.0434325 |
| 21 | 0.0544380 | 0.0249299 | 0.0793679 | 0.0898340 | 0.0632713 |
| 22 | 0.0860472 | 0.0349790 | 0.121026 | 0.141996 | 0.0887758 |
| 23 | 0.108424 | 0.0473505 | 0.155774 | 0.178922 | 0.120174 |
| 24 | 0.109624 | 0.0564823 | 0.166106 | 0.180903 | 0.143350 |
| 25 | 0.0939636 | 0.0594188 | 0.153382 | 0.155059 | 0.150803 |
| 26 | 0.0682472 | 0.0542051 | 0.122452 | 0.112622 | 0.137571 |
| 27 | 0.0397189 | 0.0398622 | 0.0795811 | 0.0655444 | 0.101169 |
| 28 | 0.0183842 | 0.0234909 | 0.0418751 | 0.0303377 | 0.0596192 |
| 29 | 0.00536205 | 0.00893675 | 0.0142988 | 0.00884850 | 0.0226812 |
| Outcome probabilities | 0.605984 | 0.394016 |  |  |  |
| Mean duration |  |  | 23.9151 | 23.8664 | 23.9899 |
| Standard deviation |  |  | 2.33806 | 2.00364 | 2.77314 |

Let $R_m(l,u)$ be the discrete hypercube $R_m(l,u) = \{l,...,u\}^m$ in $m$ dimensions. Define the stopping times

$$N_1 = N_1(m,s,l,u) = \inf\{n \geq ml : X^{(n)} \in R_m(l,u)\},$$

$$N_2 = N_2(m,s,l,u) = \inf\{n \geq u+1 : \max_j \{X^{(n)}\} = u+1\},$$

and the *band or bump stopping time*,

$$N = N(m,s,l,m) = N_1 \wedge N_2.$$

Below we obtain explicit expressions for the joint probabilities $P[N = n, band] = P[N = N_1 = n]$ and $P[N = n, bump] = P[N = N_2 = n]$ in terms of general rectangular event probabilities $P_{m'}[X^{(n')} \in R']$ for suitable choices of $m', n',$ and $R'$. The generic $X^{(n')}$ has the multiple hypergeometric distribution, $X^{(n')} \sim H_{m'}(n'; s,...,s)$, with probability function $P[X^{(n')} = x'] = \binom{s}{x'_1} \cdots \binom{s}{x'_{m'}} / \binom{t}{n'}$ for $x'$ in the support $\mathcal{H}_{m'}(n'; s,...,s) = \{x' \in \{0,...,s\}^{m'} : x'_1 + \cdots + x'_{m'} = n'\}$. We begin by assuming $0 < l < u < s$, after which we cover the boundary cases $l = u = 0$, $0 = l < u < s$, $u = s$, and $0 < l = u - s$.

−2−

It is clear that a band cannot occur before $n = ml$ draws and must occur no later than $n_{max} = l + (m-1)u$ draws (by the pidgeonhole principle). Consider the tally vector $X^{(n-1)}$ just prior to the event $[N = n, band]$. All tallies must be no larger than $u$; exactly one tally, say the $j^{th}$, must equal $l-1$, for which there are $m$ possibilities; the remaining tallies must be at least $l$; and the $n^{th}$ card must then be of type $j$, for which there are $s-(l-1)$ possibilities among the $t-(n-1)$ remaining cards. Let $x'$ denote an arbitrary vector in the discrete simplex $\mathcal{H}_{m-1}(n-l\,;s,...,s)$ and let

$$X' \sim H_{m-1}(n-l\,;s,...,s).\qquad(1)$$

The probabilities that $X' = (X_1^{(n-1)},...,\hat{X}_j^{(n-1)},...,X_m^{(n-1)}) \in R' = R_{m-1}(l,u)$ don't depend on $j$, so we have

$$P[N=n, band] = m\binom{s+1-l}{t+1-n}\sum_{x'\in R'}\prod_{j'=1}^{m-1}\binom{s}{l-1}\binom{s}{x'_{j'}}\bigg/\binom{t}{n-1}$$

$$= \left\{ m\binom{s+1-l}{t+1-n}\frac{\binom{s}{l-1}}{\binom{t}{n-1}} \right\} \binom{(m-1)s}{n-l} P_{m-1,n-l}[X' \in R'].\qquad(2)$$

But the leading term of (2) in braces is equal to

$$m\binom{s+1-l}{t+1-n}\frac{s!}{(l-1)!(s+1-l)!}\bigg/\frac{t!}{(n-1)!(t+1-n)!} = \frac{ms(s-1)!}{(l-1)!(s-l)!}\bigg/\frac{t!}{(n-1)!(t-n)!} = \binom{s-1}{l-1}\bigg/\binom{t-1}{n-1}$$

since $t = ms$, so we have

$$P[N=n, band] = \left\{\binom{s-1}{l-1}\binom{t-s}{n-l}\bigg/\binom{t-1}{n-1}\right\} P_{m-1,n-l}[X' \in R'].\qquad(3)$$

We recognize the leading term of (3) in braces as a univariate hypergeometric point probability,

$$P[Y_{n-1} = l-1] = \binom{s-1}{l-1}\binom{t-s}{n-l}\bigg/\binom{t-1}{n-1} = \binom{n-1}{l-1}\binom{t-n}{s-l}\bigg/\binom{t-1}{s-1},\qquad(4)$$

where $(Y_{n-1}, n-1-Y_{n-1}) \sim H_2(n-1;s-1,t-s)$ or, equivalently, $(Y_{n-1}, s-1-Y_{n-1}) \sim H_2(s-1;n-1,t-n)$. Thus we may write

$$P[N=n, band] = P[Y_{n-1} = l-1] P_{m-1,n-\ell}[X' \in R_{m-1}(l,u)]\qquad(5)$$

with $X'$ distributed as in (1). The latter expression in (4) together with (5) yields expression (1.2) of Levin (2024) in the case $m=13$, $s=4$, $l=1$, and $u=3$.

The hypercube event probability $P_{m-1,n-\ell}[X' \in R_{m-1}(l,u)]$ equals zero for $n-l$ outside the set $\{(m-1)l,...,(m-1)u\}$, in agreement with the range for $n$ stated above, $ml \le n \le l+(m-1)u = n_{max}$. As $\mathcal{H}_{m-1}(n-l\,;s,...,s)$ is non-empty for $n$ in this range, the hypercube event probability is also positive.

–3–

This range also ensures $P[Y_{n-1} = l-1] > 0$ since for that we would need $l-1 \leq \min(s-1, n-1)$, which is true since $l \leq s$ and $l < ml < mu \leq n$]; and $l-1 \geq \max\{0, y_{\min}\}$ with $y_{\min} = (s-1)+(n-1)-(t-1) = n-1-(m-1)s$, which is also true since $n-1-(m-1)s \leq l+(m-1)u-1-(m-1)s = l-1-(m-1)(s-u) \leq l-1$. It follows that $P[N = n, band] > 0$ for each $ml \leq n \leq n_{\max}$. Summing (5) over this range yields the marginal probability that the game terminates with a band,

$$P[band] = \sum_{n=ml}^{l+(m-1)u} P[Y_{n-1} = l-1] P_{m-1,n-l}[X' \in R_{m-1}(l,u)], \qquad (6)$$

with $X'$ and $Y_{n-1}$ distributed as in (1) and (4), respectively.

Turning now to $P[N = n, bump]$, it is clear that the only observable bump events are $[\max_j \{X_j^{(n)}\} = u+1]$ and that a bump cannot occur before $n = u+1$ and must occur no later than $n_{\max} = l+(m-1)u$ just as before, this time because if $N > l-1+(m-1)u$, the next card must result in either a band or a bump. Consider the tally vector $X^{(n-1)}$ just prior to stopping with a bump. All tallies must be no larger than $u$; at least one must be no larger than $l-1$ (else a band has occurred); and at least one must equal $u$. For any $n$ and any tally-vector $X$, let $k(X) = \sum_{j=1}^{m} I[X_j = u]$ denote the number of maximum within-quota tallies, let $k'(X) = \sum_{j=1}^{m} I[X_j < l]$ denote the number of below-quota tallies, and let $k''(X) = \sum_{j=1}^{m} I[l \leq X_j < u]$ denote the remaining within-quota tallies with $k(X) + k'(X) + k''(X) = m$. Then just prior to $[N = n, bump]$ we must have $k = k(X^{(n-1)}) \geq 1$ and $k' = k'(X^{(n-1)}) \geq 1$, and with $k'' = m-k-k'$, there are $\binom{13}{k,k',k''}$ ways to choose the corresponding subsets of tallies. Furthermore, the rank of the $n^{th}$ card drawn must be among the $k$ tallies equal to $u$ to complete the bump, for which there are $k(s-u)$ possibilities among the $t-(n-1)$ remaining cards. Now let $(x', x'')$ denote an arbitrary vector in $\mathcal{H}_{m-k}(n_k; s,...,s)$ and let

$$(X', X'') \sim H_{m-k}(n_k; s,...,s), \text{ where} \qquad (7)$$

$$n_k = n-1-ku. \qquad (8)$$

The probability that the tallies not equal to $u$ after $n-1$ draws, which sum to $n_k$, all lie in $R(k', k'') = R_{k'}(0, l-1) \times R_{k''}(l, u-1)$ doesn't depend on the particular subsets of tallies, it depends only on the *configuration* $(k, k', k'')$. Thus we have

$P[N = n$ with a bump from configuration $(k, k', k'')$ after $n-1$ draws$]$

$$= \binom{m}{k,k',k''}\left(\frac{k(s-u)}{t+1-n}\right)\binom{s}{u}^k \sum_{(x',x'')\in R(k',k'')} \prod_{j'}\binom{s}{x'_{j'}}\prod_{j''}\binom{s}{x''_{j''}} \bigg/ \binom{t}{n-1} \tag{9}$$

$$= \left\{\binom{m}{k,k',k''}\left(\frac{k(s-u)}{t+1-n}\right) \bigg/ \binom{t}{n-1}\right\}\binom{s}{u}^k \binom{t-ks}{n_k} P_{n_k}[(X',X'') \in R(k',k'')].$$

The leading term of (9) in braces is equal to

$$\left(\frac{m!}{k!k'!k''!}\right)\left(\frac{k(s-u)}{t+1-n}\right) \bigg/ \left\{\frac{t!}{(n-1)!(t+1-n)!}\right\} = \left(\frac{m!}{k!(m-k)!}\frac{(m-k)!}{k'!k''!}\right)\left(\frac{k(s-u)}{n}\right) \bigg/ \left\{\frac{t!}{n!(t-n)!}\right\}$$

$$= \binom{m}{k}\frac{k(s-u)}{n\binom{t}{n}}\binom{m-k}{k''}.$$

Thus, summing over all configurations $(k, k', k'')$, we find

$$P[N = n, bump] = \sum_k \left\{\binom{m}{k}\frac{k(s-u)}{n\binom{t}{n}}\right\}\binom{s}{u}^k\binom{t-ks}{n_k}\sum_{k''}\binom{m-k}{k''}P_{m-k,n_k}[(X',X'')\in R(k',k'')] \tag{10}$$

with $(X', X'')$ distributed as in (7). We specify the ranges of $n$, $k$, and $k''$ that are necessary to yield positive summands in (10) as follows.

As mentioned above, $n$ must satisfy

$$u+1 \le n \le n_{\max} = l + (m-1)u \tag{11}$$

For the outer summation, $k$ must satisfy three constraints. First, $k \ge 1$. Second, no matter what values $k'$ or $k''$ may take, we must have $n_k \le k'(l-1) + k''(u-1)$, else no $(X'^{(n-1)}, X''^{(n-1)})$ could satisfy $\sum_{j'} X'^{(n-1)}_j + \sum_{j''} X''^{(n-1)}_j = n_k$. But since $k' = m - k - k''$ must be no less than 1, we have $k'' \le m - k - 1$, so

$$n - 1 - ku = n_k \le (m-k-k'')(l-1) + k''(u-1) = (m-k)(l-1) + k''(u-l) \le (m-k)(l-1) + (m-k-1)(u-l),$$

which implies $k \ge n - \{l + (m-1)(u-1)\} = n - \{n_{\max} - (m-1)\}$. Third, we must have $n - 1 \ge ku$, so that $k \le \lfloor (n-1)/u \rfloor$. Therefore, we may restrict the range of $k$ in the outer summation of (10) to

$$\max\{1, n - \{l + (m-1)(u-1)\}\} \le k \le \lfloor (n-1)/u \rfloor. \tag{12}$$

For the inner summation, for any given $n$ and $k$ satisfying (11) and (12), respectively, $k''$ must satisfy two constraints in addition to $k'' \ge 0$ and $k'' \le m - k - 1$ as already indicated. Note that we must have $n-1 \leq ku + k'(l-1) + k''(u-1)$, else no $X^{(n-1)}$ could satisfy $\sum_{j=1}^{m} X_j^{(n-1)} = n-1$. This requires that $n_k = n-1-ku \leq (m-k-k'')(l-1) + k''(u-1) = (m-k)(l-1) + k''(u-l)$. Also, because $x''_{j''} \geq l$ for $(x',x'') \in R(k',k'')$, we must have $n_k \geq \sum_{j''} X_{j''}^{(n-1)} \geq k''l$, so that $k'' \leq \lfloor n_k/l \rfloor$. Therefore, we may restrict the range of $k''$ in the inner summation of (10) to

$$\max\left\{0, \left\lceil \frac{n_k - (m-k)(l-1)}{u-l} \right\rceil \right\} \leq k'' \leq \min\{\lfloor n_k/l \rfloor, m-k-1\}. \tag{13}$$

In the case $m=13$, $s=4$, $l=1$, $u=3$ of the Book or Band game, $R(k',k'') = R_{k'}(0,0) \times R_{k''}(1,2) = \{0\}^{k'} \times \{1,2\}^{k''}$, $\prod_{j'} \binom{4}{0} = 1$ in the first line of (9), and we may renormalize the sum of terms over only $k''$ by $\binom{4k''}{n_k}$ instead of over both $k'$ and $k''$ by $\binom{t-ks}{n_k}$, in which case (10) becomes

$$P[N=n, bump] = \sum_{k} \binom{13}{k} \left[ \frac{k4^k}{n\binom{52}{n}} \right] \sum_{k''} \binom{13-k}{k''} \binom{4k''}{n_k} P_{k'',n_k}[X'' \in \{1,2\}^{k''}],$$

which is equivalent to equation (1.3) of Levin (2024). From (12), $k$ ranges from $\max\{1, n-25\}$ to $\lfloor (n-1)/3 \rfloor$ and from (13), $k''$ ranges from $\lceil n_k/2 \rceil$ to $\min\{n_k, 12-k\}$. When $k''=0$, we have $n_k = 0$ and (1.3) simplifies to the sum over $k$ of terms $\binom{13}{k} \frac{k4^k}{n\binom{52}{n}}$. When $k''=1$, we have $n_k = 1$ or $2$ and (1.3) simplifies to a sum over $k$ of terms $\binom{13}{k,12-k,1} \frac{k4^k \binom{4}{n_k}}{n\binom{52}{n}} = \frac{13\binom{12}{k}k4^k\binom{4}{n_k}}{n\binom{52}{n}}$ since $P_{1,n_k}[X'' \in \{1,2\}] = 1$.

**Remark 1.** We derived the constraints in (13) on $k''$ as necessary for the summands in (10) to be positive. This leaves open the non-obvious question of whether (11) and (12) are also sufficient to guarantee that (13) produces a *non-vacuous* range for $k''$ with positive summands in (10). We show this is the case in the Appendix. □

We conclude the derivation by providing expressions for $P[N=n, band]$ and $P[N=n, bump]$ in the boundary cases of $l$ and $u$ not previously considered.

(a) $l = u = 0$. Only a bump can occur, at $N=1$ with probability 1.

(b) $0 = l < u$. Only a band can occur, at $N=1$ with probability 1.

(c) $u = s$. Only a band can occur, at time $N_l = \inf\{n \geq ml : \min_j\{X_j^{(n)}\} = l\}$. This is a coupon collector's problem with a uniform but finite population of coupons and a minimum demand of $l$ coupons per type. Because $P[N_l \leq n] = P[X^{(n)} \in R_m(l,s)]$, we have, with $X^{(n)} \sim H_m(n;s,...,s)$,

$$P[N = n, band] = P[N_l = n] = P[X^{(n)} \in R_m(l,s)] - P[X^{(n-1)} \in R_m(l,s)].$$

(d) $0 < l = u < s$. A band can only occur at time $n_{max} = l + (m-1)u = mu$, though bumps can occur sooner. Let $N_u = \inf\{n \geq u+1 : \max_j\{X_j^{(n)}\} > u\}$, for which $P[N_u > n] = P[X^{(n)} \in R_m(0,u)]$. For any $n \leq mu$, $X^{(n-1)} \in [N_u > n-1]$ means no bump has occurred by time $n-1$, so $X_j^{(n)} \leq u$ for each $j$, and if $X^{(n)} \notin [N_u > n]$ as well, then $X_j^{(n)} = u+1$ for some $j$ for the first time. Thus

$$P[N = n, bump] = P[X^{(n-1)} \in R_m(0,u)] - P[X^{(n)} \in R_m(0,u)]$$

(and $P[N = n, band] = 0$ for $n < mu$). For $n = n_{max}$, $X^{(n)} \in [N_u > n-1]$ implies $X^{(n-1)}$ is some permutation of $(u,...,u,u-1)$, so the next card must result in either a band or a bump. Thus $P[N_u > n_{max} - 1] = P[N = n_{max}] = P[N = n_{max}, band] + P[N = n_{max}, bump]$, whence

$$P[N = n_{max}, band] = P[N_u > n_{max} - 1] - P[N = n_{max}, bump] = P[X^{(n_{max})} \in R_m(0,u)].$$

**Remark 2.** A non-negative sequence $\{p_n : n = 0,...,t\}$ is called (discrete) *log-concave* if $p_{n+1} p_{n-1} \leq p_n^2$ for $0 < n < t$ and the support $\{n : p_n > 0\}$ consists of consecutive integers. Levin (2024) proves that for any $m$-vector $s = (s_1,...,s_m)$ of positive integers, any integer vectors $l = (l_1,...,l_m)$ and $u = (u_1,...,u_m)$ with $0 \leq l_j \leq u_j \leq s_j$, and multiple hypergeometric $X^{(n)} \sim H_m(n; s_1,...,s_m)$, the sequence $\{P[X^{(n)} \in R_m(l,u)] : n = 0,...,t\}$ is log-concave in $n$, where $R_m(l,u) = \{l_1,...,u_1\} \times \cdots \times \{l_m,...,u_m\}$ are arbitrary $m$-dimensional rectangles. It follows that the sequence $\{P[N = n, band] : n = ml,...,n_{max}\}$ in (5) is log-concave, since $P[Y_{n-1} = l-1] = \binom{s-1}{l-1}\binom{t-s}{n-l}/\binom{t-1}{n-1}$ is log-concave in $n$ and the product of two log-concave sequences is log-concave (Saumard and Wellner, 2014). It is not clear how to give a similar argument for $\{P[N = n, bump] : n = ml,...,n_{max}\}$ in (10), though we conjecture that (10) is also log-concave based on extensive numerical evidence. □

**Appendix**

Here we show that constraints (11) and (12) on $n$ and $k$, respectively, imply non-vacuous ranges for $k''$ in (13). The approach will be to identify successive domains of $m$, $l$, $u$, $n$, and $k$ leading to non-vacuous ranges in (13) until no possible combinations remain that could render (13) vacuous. We do this in four steps. Let

–7–

$$k_L = k_L(n) = n - \{l + (m-1)(u-1)\},$$

so we may rewrite (12) as

$$\max\{1, k_L\} \leq k \leq \lfloor (n-1)/u \rfloor; \tag{A.1}$$

and let

$$k''_L = k''_L(n, k) = \frac{n_k - (m-k)(l-1)}{u - l},$$

so that we may rewrite (13) as

$$\max\{0, \lceil k''_L \rceil\} \leq k'' \leq \min\{\lfloor n_k / l \rfloor, m - k - 1\}. \tag{A.2}$$

<u>Step 1</u>. (A.2) is never vacuous if $k''_L \leq 0$, for in that case its left-hand side is 0 while its right-hand side is never less than zero, since $m - k - 1 = k' + k'' - 1 \geq 0$. So we assume henceforth that $k''_L > 0$, which is equivalent to

$$k < \frac{n - 1 - m(l-1)}{u - l + 1} \tag{A.3}$$

or equivalently,

$$k \leq \left\lceil \frac{n - 1 - m(l-1)}{u - l + 1} \right\rceil - 1. \tag{A.4}$$

In particular, if $k = m - 1$, then $n - 1 \leq l - 1 + (m-1)u$, since $k'$ must equal 1 with $k'' = 0$, so that $n_k = n - 1 - (m-1)u \leq l - 1$, in which case $k''_L \leq 0$. So we may assume henceforth that

$$k \leq m - 2. \tag{A.5}$$

<u>Step 2</u>. We claim that

$$m - k - 1 \geq \lceil k''_L \rceil. \tag{A.6}$$

This is because the following equivalent statements show that

$$k \geq k_L \text{ holds if and only if } m - k - 1 \geq k''_L.$$

$$m - k - 1 \geq k''_L$$

$$(m - k - 1)(u - l) \geq n_k - (m - k)(l - 1)$$

$$(m - 1)(u - l) - k(u - l) \geq n - 1 - ku - (m - k)(l - 1)$$

$$(m - 1)(u - l) \geq n - 1 - (m - k)(l - 1) - ku + k(u - l) = n - 1 - m(l - 1) - k, \text{ or}$$

$$k \geq n - 1 - m(l - 1) - (m - 1)(u - l)$$

$$= n - 1 - m(l - 1) - (m - 1)\{(u - 1) - (l - 1)\}$$

$$= n - 1 - m(l - 1) + (m - 1)(l - 1) - (m - 1)(u - 1)$$

$$= n - \{l + (m - 1)(u - 1)\} = k_L.$$

Thus $m - k - 1 \geq k''_L$ since $k \geq k_L$ under (A.1) and because $m - k - 1$ is an integer, $m - k - 1 \geq \lceil k''_L \rceil$. □

A consequence of (A.6) is that whenever $\lfloor n_k/l \rfloor \geq m-k-1$, the right-hand side of (A.2) equals $m-k-1 \geq \lceil k_L'' \rceil$, which is the left-hand side of (A.2), so (A.2) is non-vacuous in such cases. So we assume henceforth that

$$\lfloor n_k/l \rfloor < m-k-1. \tag{A.7}$$

(A.7) is equivalent to

$$n_k/l < m-k-1 \tag{A.8}$$

because (A.8) implies (A.7) and conversely, if (A.7) holds but not (A.8), then $n_k/l \geq m-k-1 \in \mathbb{Z}$ implies $\lfloor n_k/l \rfloor \geq m-k-1$, contradiction. Rearranging terms shows that (A.8) is equivalent to

$$k > \frac{n-1-(m-1)l}{u-l} \tag{A.9}$$

or equivalently,

$$k \geq \left\lfloor \frac{n-1-(m-1)l}{u-l} \right\rfloor + 1. \tag{A.10}$$

Under (A.7), the right-hand side of (A.2) equals $\lfloor n_k/l \rfloor$. In particular, when $l=1$, we have $u \geq 2$ and $k_L'' = n_k/(u-1) = (n-1-ku)/(u-1) \leq n-1-ku = \lfloor n_k/l \rfloor$, in which case $\lceil k_L'' \rceil \leq \lfloor n_k/l \rfloor$ and (A.2) is non-vacuous. So we assume henceforth that

$$l \geq 2. \tag{A.11}$$

Step 3. We note that not every $n \leq n_{\max} = l+(m-1)u$ admits of $k$ satisfying (A.4) and (A.10), but that only means such $k$ satisfying (A.1) have already been shown to imply non-vacuous (A.2). In fact, it is clear from (A.3) and (A.9) that if $\frac{n-1-(m-1)l}{u-l} \geq \frac{n-1-m(l-1)}{u-l+1}$ then no $k$ can satisfy (A.4) and (A.10). This is because for any $x$, $\lfloor x \rfloor + 1 > x$ and $\lceil x \rceil - 1 < x$ and so

$$\left\lfloor \frac{n-1-(m-1)l}{u-l} \right\rfloor + 1 > \frac{n-1-(m-1)l}{u-l} \geq \frac{n-1-m(l-1)}{u-l+1} > \left\lceil \frac{n-1-m(l-1)}{u-l+1} \right\rceil - 1.$$

But $\frac{n-1-(m-1)l}{u-l} \geq \frac{n-1-m(l-1)}{u-l+1}$ holds if and only if $n-1 \geq (m-l)u + l(l-1)$, so we need only consider values of $n-1 < (m-l)u + l(l-1)$, i.e.,

$$n \leq n_{HI} = (m-l)u + l(l-1) = n_{\max} - (u-l)(l-1) - l. \tag{A.12}$$

Henceforth we assume (A.12). Smaller values of $n$ may also not admit of $k$ satisfying (A.4) and (A.10), but no $n > n_{HI}$ admits of such $k$.

Step 4.  We now claim that under assumptions (A.1), (A.5), (A.11), and (A.12), for any $n$ which admits of $k$ satisfying (A.4) and (A.10), we have

$$\frac{n_k}{l} - k_L'' \geq 1. \qquad (A.13)$$

(A.13) implies (A.2) is non-vacuous, because $n_k/l = \lfloor n_k/l \rfloor + \varepsilon_1$ and $k_L'' = \lceil k_L'' \rceil - \varepsilon_2$, with $0 \leq \varepsilon_1, \varepsilon_2 < 1$, so $\lfloor n_k/l \rfloor - \lceil k_L'' \rceil = (n_k/l) - k_L'' \geq 1 - (\varepsilon_1 - \varepsilon_2) \geq 0$. To prove the claim, first note that (A.13) holds if and only if the following equivalent statements hold.

$$\frac{n_k}{l} - \frac{n_k - (m-k)(l-1)}{u-l} \geq 1$$

$$n_k(u-l) - l\{n_k - (m-k)(l-1)\} \geq l(u-l)$$

$$n_k(u-2l) + (m-k)l(l-1) \geq l(u-l)$$

$$(n-1-ku)(u-2l) + (m-k)l(l-1) \geq l(u-l)$$

$$k\{u(u-2l) + l(l-1)\} \leq (n-1)(u-2l) + ml(l-1) - l(u-l).$$

This yields the necessary and sufficient condition for (A.13), namely,

$$k\{(u-l)^2 - l\} \leq (n-1)(u-2l) + l\{m(l-1) - (u-l)\}. \qquad (A.14)$$

So we must demonstrate that any $k$ satisfying (A.4) and (A.10) also satisfies (A.14). There are three cases to consider: (a) $(u-l)^2 \geq l$, $u \geq 2l$; (b) $(u-l)^2 \geq l$, $u < 2l$; and (c) $(u-l)^2 < l$, $u < 2l$. These are the only cases since $(u-l)^2 < l$ implies $u = u - l + l \leq (u-l)^2 + l < 2l$.

In case (a), for any given $k$ satisfying (A.4), the right-hand side of (A.14) is least when $n$ is the least possible value such that the right-hand side of (A.4) equals $k$, i.e., when $\left\lceil \frac{n-1-m(l-1)}{u-l+1} \right\rceil - 1 = k$, or equivalently, when $n$ is one more than the greatest possible value such that $\left\lceil \frac{n-1-m(l-1)}{u-l+1} \right\rceil - 1 = k - 1$, which occurs when $n - 1 = m(l-1) + k(u-l+1)$. Thus (A.14) is least for $n - 1 = 1 + m(l-1) + k(u-l+1)$, in which case (A.14) becomes

$$k\{(u-l)^2 - l\} \leq \{1 + m(l-1) + k(u-l+1)\}(u-2l) + l\{m(l-1) - (u-l)\}. \qquad (A.15)$$

Collecting terms in $k$ on the left, the coefficient of $k$ is $(u-l)^2 - l - (u-l+1)(u-2l)$

$$= (u-l)^2 - l - (u-l+1)(u-l-l) = (u-l)^2 - l - (u-l)^2 + l(u-l) - (u-l) + l = (u-l)(l-1).$$

On the right, we have $\{1 + m(l-1)\}(u-2l) + l\{m(l-1) - (u-l)\}$

–10–

$$= m\{(l-1)(u-2l)+l(l-1)\}+(u-2l)-l(u-l) = m\{(l-1)(u-l-1)+l(l-1)\}+(u-l-1)-l(u-l)$$
$$= m(u-l)(l-1)-(u-l)(l-1)-l.$$

Thus, with $l > 1$ under (A.11), the condition on $k$ in (A.14) becomes $k \le m-1-\dfrac{l}{(u-l)(l-1)}$ and any $k$ satisfying that condition would also satisfy (A.14) for any larger value of $n$ yielding the same value of $k$ on the right-hand side of (A.4). But since $u-l \ge l$ in case (a), we have $0 < l/\{(u-l)(l-1)\} < 1$, so the condition becomes $k \le m-2$, which holds by (A.5).

In case (b), $(u-l)^2 \ge l$, $u < 2l$, for any given $k$ satisfying (A.4), the right-hand side of (A.14) is least when $n$ is the greatest possible value such that the right-hand side of (A.4) equals $k$. This occurs when $n-1$ is $u-l$ more than the previously identified least value. Then the coefficient of $k$ on the left-hand side (A.15) remains the same while the right-hand side has the additional term $(u-l)(u-2l)$, so the condition on $k$ in (A.14) becomes $k \le m-1-\left\{\dfrac{l+(u-l)(2l-u)}{(u-l)(l-1)}\right\}$. Now the term in braces is no greater than 2, since that is true if and only if $l \le 2(u-l)(l-1)-(u-l)(2l-u) = (u-l)(u-2)$, while $l \ge 2$ implies $(u-l)(u-2) \ge (u-l)^2 \ge l$. Therefore it will suffice to show in this case that $k \le m-3$ for any $k$ satisfying (A.4). But under (A.12), the greatest possible value on the right-hand side of (A.4) occurs at $n = n_{HI}$ with value $m-l-1$. This is because at $n = n_{HI} = u(m-l)+l(l-1)$,

$$\dfrac{n_{HI}-1-m(l-1)}{u-l+1} = \dfrac{u(m-l)+l(l-1)-1-m(l-1)}{u-l+1} = \dfrac{u(m-l)-(l-1)(m-l)-1}{u-l+1} = m-l-\dfrac{1}{u-l+1},$$

so $\left\lceil\dfrac{n_{HI}-1-m(l-1)}{u-l+1}\right\rceil - 1 = m-l-1$. Thus, since $l \ge 2$, $k$ in (A.4) can be no greater than $m-3$.

In case (c), $(u-l)^2 < l$, $u < 2l$, necessary and sufficient condition (A.14) becomes

$$k \ge \dfrac{(n-1)(2l-u)-l\{m(l-1)-(u-l)\}}{l-(u-l)^2}.$$

Letting $d = u-l \ge 1$ with $2l-u = l-d > 0$ and $l-(u-l)^2 = l-d^2 > 0$, the condition becomes

$$k \ge \dfrac{(n-1)(l-d)-l\{m(l-1)-d\}}{l-d^2}. \quad (A.16)$$

But any $k$ satisfying (A.10) exceeds $\{n-1-l(m-1)\}/(u-l) = \{n-1-l(m-1)\}/d$, so it suffices to show for such $k$ that

$$\dfrac{n-1-l(m-1)}{d} \ge \dfrac{(n-1)(l-d)-l\{m(l-1)-d\}}{l-d^2}.$$

–11–

Rearranging terms, we are to show the equivalent inequalities

$$(n-1)\{\frac{l-d}{l-d^2} - \frac{1}{d}\} \geq \frac{ml(l-1)(n-1) - ld}{l-d^2} - \frac{l(m-1))\}}{d}$$

$$(n-1)\{d(l-d) - l + d^2\} \leq d\{ml(l-1) - ld\} - (m-1)l(l-d^2)$$

$$(n-1)(d-1) \leq d\{m(l-1) - d\} - (m-1)(l-d^2)$$

$$= m(d^2 + dl - d - l) + l - 2d^2$$

$$= m(d+l)(d-1) + l - 2d^2 = mu(d-1) + l - 2d^2.$$

When $d=1$, this is obvious as $l \geq 2$ by (A.11). When $d > 1$ we are to show that

$$n - 1 \leq mu + \frac{l - 2d^2}{d-1}.$$

Now by (A.12), $n \leq n_{HI} = (m-l)u + l(l-1) = mu - l(u-l+1) = mu - l(d+1)$, so it suffices to show that

$$\frac{l - 2d^2}{d-1} + l(d+1) + 1 \geq 0,$$

or equivalently,

$$l - 2d^2 + l(d^2 - 1) + d - 1 = d^2(l-2) + d - 1 \geq 0,$$

which is true by (A.12) and $d \geq 1$.

This concludes the proof of (A.13) and thus (A.2) is non-vacuous for all $n$ satisfying (11) and $k$ satisfying (A.1)=(12).  □